\documentclass[12pt]{article}

\usepackage[english]{babel}

\usepackage[a4paper,top=2cm,bottom=2cm,left=3cm,right=3cm,marginparwidth=1.75cm]{geometry}

\usepackage{amsmath, authblk}
\usepackage{graphicx}
\usepackage[colorlinks=false]{hyperref}

\usepackage{blindtext}
\usepackage{parskip}
\usepackage{listings}
\usepackage{amsfonts}
\usepackage{amsthm}
\usepackage{amssymb}
\usepackage{fullpage}
\usepackage{mathtools}
\usepackage{enumerate}
\usepackage{minted}
\usepackage{tikz}
\usetikzlibrary{matrix}
\usepackage{lscape}
\usepackage{longtable}

\theoremstyle{plain}
\newtheorem{theorem}{Theorem}[section]

\newtheorem{lemma}[theorem]{Lemma}

\theoremstyle{definition}
\newtheorem{definition}[theorem]{Definition}
\usepackage{graphicx}

\newtheorem{teo}{Theorem}  

\usepackage{csquotes}
\title{On the convergence of boundary points for hyperbolic inner functions}

\author[1]{Anna Jov\'e\thanks{Supported by the Spanish government grant FPI PRE2021-097372 and PID2023-147252NBI00.}}
\author[2]{Mateo Menc\'ia\thanks{Supported by Fundación María Cristina Masaveu Peterson. Corresponding author. \url{mateomenciarodriguez@gmail.com}}}
\affil[1]{\small Departament de Matemàtiques i Informàtica, Universitat de Barcelona, Barcelona, Spain}
\affil[2]{\small Faculty of Sciences, University of Oviedo, Oviedo, Spain}
\begin{document}
\maketitle

\begin{abstract}
Given a hyperbolic inner function \( f\colon \mathbb{D}\to\mathbb{D} \) with Denjoy--Wolff point \( p\in\partial\mathbb{D} \), it is well known that almost every point \( \xi\in \partial\mathbb{D} \) converges to \( p \) under iteration of the radial extension \( f^*\colon\partial\mathbb{D}\to\partial\mathbb{D} \). We provide explicit bounds for the rate of this convergence in terms of the angular derivative, holding almost surely. Our results also cover the case where the Denjoy--Wolff point is a singularity.
\end{abstract}

\section{Introduction}

Let $\mathbb{D}$ be the unit disk,  let $f \colon \mathbb{D} \to \mathbb{D}$ be 
a holomorphic self map of $\mathbb{D}$, and consider the dynamical system given by its iterates $\{ f^n \}_{n\in \mathbb{N}}$.
The Denjoy-Wolff Theorem states that if $f$ is not conjugate to a rotation, then there exists a point $p \in \overline{\mathbb{D}}$ such that the orbit of every point in $\mathbb{D}$ converges to $p$, that is, $f^n(z) \to p$ as $n \to \infty$ for every $z \in \mathbb{D}$. We say that $p$ is the {\em Denjoy-Wolff point} of $f$.

In the case where $f$ is an {\em inner function}, that is, for $\lambda$-almost every $\zeta\in\partial \mathbb{D}$ the radial limit \[ f^*(\zeta)=\lim_{r\to1} f(r\zeta)\] belongs to $\partial\mathbb{D}$ (where $\lambda$ denotes the normalized Lebesgue measure on $\partial\mathbb{D}$), one can consider the dynamical system defined in the unit circle given by the radial extension
$$
f^{*} \colon \partial\mathbb{D} \longrightarrow \partial\mathbb{D}.
$$
Assuming that the Denjoy-Wolff point $p$ lies in $\partial\mathbb{D}$, a natural question to ask is whether points in the unit circle also converge to the Denjoy-Wolff point under the iteration of $f^*$.
In the seminal work of Aaronson, Doreing and Mañé \cite{Aaronson78, mane_dynamics_1991},  this question is answered by means of a complete characterization in terms of infinite sums. More precisely, $\lambda$-almost every point on $\partial\mathbb{D}$ converges to the Denjoy-Wolff point $p\in\partial \mathbb{D}$ if and only if \[\sum_{n\geq 0} 1-|f^n(0|<\infty.\]

It is well-known that {\em hyperbolic} inner functions (i.e. for which the angular derivative $f'(p)\in (0,1)$, see Section \ref{sect-hyperbolic-inner}) always satisfy the condition above. Going one step further, one may ask  at which
rate do these orbits approach the Denjoy–Wolff point.

If the Denjoy-Wolff point $p\in\partial \mathbb{D}$ is not a singularity (i.e. $f$ extends as a holomorphic map around $p$), then $p$ is an attracting fixed point for $f$. Using Koenigs' linearizing coordinates, we see that  $f$ behaves like the map $z \mapsto f'(p) z$ near $p$, showing geometric convergence to the fixed point for points on $\partial \mathbb{D}$ in a neighborhood of  $p$.

The case when the Denjoy–Wolff point is a singularity is much more complicated, due to the lack of normal forms. In this paper, we establish an explicit rate of convergence for hyperbolic inner functions, covering the case when the Denjoy–Wolff point is a singularity. Namely, if we denote $A(p;a,b) \coloneqq \{ z \in \mathbb{C}\colon a <|z-p| < b\}$, we prove the following.   \\
\begin{teo}\label{A}
        Let $f \colon \mathbb{D} \to \mathbb{D}$ be a hyperbolic inner function with Denjoy-Wolff point $p \in \partial{\mathbb{D}}$. Let $\alpha = f'(p) \in (0,1)$. Then for all $0<\varepsilon<1$, we have $$(f^*)^n(\zeta) \in A(p;\alpha^{(1+\varepsilon)n} , \alpha^{(1-\varepsilon)n}) \cap \partial\mathbb{D}$$ for $n$ large enough and $\lambda$-almost every $\zeta \in \partial\mathbb{D}$.
\end{teo}\
The proof of the theorem is  based on the concept of \textit{shrinking targets} of $\partial\mathbb{D}$ and their hitting properties (see Section \ref{sec:shrinking}), developed in \cite{benini_shrinking_2024} to analyze the recurrent behavior of compositions of inner functions fixing 0. 

Without being precise, we say that a shrinking target is a collection of  arcs of $\partial\mathbb{D}$ shrinking in length, and the problem is to determine whether orbits hit these arcs almost surely. Indeed, the arcs  $A(p;\alpha^{(1+\varepsilon)n} , \alpha^{(1-\varepsilon)n})$ in Theorem \ref{A} form a shrinking target, and we want to conclude that it is hit almost surely. Via Möbius transformations, we convert our autonomous system into a non-autonomous one fixing $0$, to be able to apply the criteria in \cite{benini_shrinking_2024} to determine that the shrinking target is hit almost surely. We make use of some facts about the rate of convergence of the orbit of $0$ to the Denjoy-Wolff point for hyperbolic inner functions.
\

{\bf Acknowledgements. } The authors gratefully acknowledge the Barcelona Introduction to Mathematical Research (BIMR) program at the Centre de Recerca Matemàtica (CRM) for providing an excellent research environment and support during the development of this work.

\section{Basic background on shrinking targets}\label{sec:shrinking}

In this section, we recall several results on shrinking targets, primarily from \cite{benini_shrinking_2024}, which will be instrumental in proving our main theorem. We begin by precisely defining a shrinking target and what it means for an inner function to hit such a target.\\

 \begin{definition}
     \begin{enumerate}
         \item[(i)] A {\em shrinking target} is a sequence $(I_n)$ of arcs of $\partial\mathbb{D}$, not necessarily nested, such that $|I_n| \to 0$ as $n \to \infty$, where $|I_n|$ denotes the length of the arc $I_n$.
         \item[(ii)]      Let $(F_n)$ be a sequence of inner functions, and let $\zeta \in \partial\mathbb{D}$. The sequence $(F_n(\zeta))$ {\em hits the target} $(I_n)$ if we have $F_n^*(\zeta) \in I_n$ infinitely often.
     \end{enumerate}

 \end{definition}

The following lemma provides a sufficient condition for a sequence of  inner functions not to hit a shrinking target almost surely, expressed solely in terms of the measure of the targets. \\

\begin{lemma}[{\cite[Theorem~A]{benini_shrinking_2024}}]\label{BEF+A}
    Let $(F_n)$ be a sequence of holomorphic self maps of $\mathbb{D}$ such that $F_n(0) = 0$ for $n \in \mathbb{N}$ and suppose that $(I_n)$ is a shrinking target in $\partial\mathbb{D}$ such that
    $$
    \sum_{n=1}^\infty |I_n| < \infty.
    $$
    Then $(F_n(\zeta))$ fails to hit $(I_n)$ for almost all $\zeta \in \partial \mathbb{D}$.
\end{lemma}\

In order to apply the previous lemma, we first need to map, via a Möbius transformation, $f^n(0)$ to $0$, and to study the length of the resulting arcs. For this purpose, we will take advantage of Lemma \ref{Mn}.\\

\begin{lemma}[{\cite[Theorem 5.1]{benini_shrinking_2024}}]\label{Mn}
    Let $f \colon \mathbb{D} \to \mathbb{D}$ be a inner function. Let 
    $$
M_n(z) = \frac{z+f^n(0)}{1+\overline{f^n(0)}z}, \text{\ for $n \in \mathbb{N}$},
$$
    and $M_0(z) = z$. Then, for every arc $J_n = (e^{i\theta_n},e^{i\varphi_n}) \subset \partial \mathbb{D}$, with $\theta_n < \varphi_n < \theta +2\pi$, if $I_n = M_n^{-1}(J_n)$, $n \in \mathbb{N}$, we have
    \begin{equation*}\label{senoarcoeng}
    2 \sin\left(\frac{1}{2}|I_n|\right) 
    = |M_n^{-1}(e^{i\varphi_n}) - M_n^{-1}(e^{i\theta_n})| 
    =\frac{(1-|f^n(0)|^2)|e^{i\varphi_n}-e^{i\theta_n}|}{|f^n(0) - e^{i\varphi_n}| |f^n(0) - e^{i\theta_n}|}.
     \end{equation*}
\end{lemma}\

\section{Hyperbolic self-maps of $\mathbb{D}$}\label{sect-hyperbolic-inner}

Let us recall the concept of Stolz angle.\\

\begin{definition}
    A \textit{Stolz angle} at $a \in \partial\mathbb{D}$ is a region
    $$
    \Delta = \{ z \in \mathbb{D}\colon |\arg(1-\overline{a}z)| < \alpha, |z-a| < p\}
    $$
    with $0 < \alpha < \frac{\pi}{2}$ and $\rho < \cos(2\alpha)$.
\end{definition}\ 

It is convenient to define limits and derivatives within Stolz angles, which are called angular limits and derivatives, respectively.\\

\begin{definition}\label{angularderivative}
    Let $f \colon \mathbb{D} \to \mathbb{C}$ be a holomorphic function. We say that $f$ has \textit{angular limit} $b \in \mathbb{C}\cup \{\infty\}$ at $a \in \partial\mathbb{D}$ if $\lim_{z\to a, z \in \Delta}f(z) = b$ for every $\Delta$ Stolz angle at $a$.

  We say that $f$ has \textit{angular derivative} $b \in \mathbb{C} \cup \{ \infty \}$ at $a \in \partial\mathbb{D}$ if the angular limit $f(a) \neq \infty$ exists and 
    $$
    \lim_{z\to a, z \in \Delta} \frac{f(z)-f(a)}{z-a} = b
    $$
    for every Stolz angle $\Delta$ at $a$. We write $b= f'(a)$.
\end{definition}\

As a consequence of the Julia-Wolff lemma {\cite[Proposition~4.13]{pommerenke_distortion_1992}}, if we have $f \colon \mathbb{D} \to \mathbb{D}$ holomorphic with Denjoy-Wolff point $p \in \mathbb{D}$, then $f'(p)$ exists,  and $0 < f'(p) \leq 1$. \\

\begin{definition}
    Let $f\colon \mathbb{D} \to \mathbb{D}$ be holomorphic, with Denjoy-Wolff point $p \in \partial\mathbb{D}$. We say that $f$ is \textit{hyperbolic} if $f'(p) \in (0,1)$.
\end{definition}\

For self-maps of $\mathbb{D}$ with Denjoy-Wolff point on $\partial \mathbb{D}$, the behavior of orbits inside $\mathbb{D}$ is well understood, as shown by Wolff lemma.\\

\begin{lemma}[Wolff lemma, {\cite[Section~4.4]{shapiro_angular_1993}}]\label{cotas}
    Let $f\colon \mathbb{D} \to \mathbb{D}$ be holomorphic, with Denjoy-Wolff point $p \in \partial\mathbb{D}$.   Consider $\eta >0$ and 
    $$
    H(p,\eta) = \{ z\in \mathbb{D} \colon |p-z|^2 < \eta (1-|z|)^2\}.
    $$
    Then, $f(H(p,\eta)) \subseteq H(p,f'(p)\eta)$.
\end{lemma}\

Moreover, it is known that for hyperbolic self maps of $\mathbb{D}$, the convergence of interior orbits to the Denjoy-Wolff point is {\em non-tangential}, that is, within a Stolz angle (see \cite{cowen_iteration_1981}). 

Finally, regarding boundary orbits of hyperbolic inner functions, note that by combining {\cite[Theorem~3.1, Corollary~4.3]{mane_dynamics_1991} and \cite[Lemma~2.6]{stallard_iteration_2008}}, it can be proved that if $f$ is a hyperbolic inner function, then $\sum_{n=0}^\infty (1-|f^n(0)|) < \infty$. By \cite[Theorem~4.1]{mane_dynamics_1991}, this means that $f^n(\zeta) \to p$ for $\lambda$-almost every $\zeta \in \mathbb{D}$. Our results prove the exact rates of convergence.

\section{Proof of Theorem \ref{A}}

In this section, we prove Theorem \ref{A}, which states that if $f \colon \mathbb{D} \to \mathbb{D}$ is a hyperbolic inner function with Denjoy-Wolff point $p \in \partial{\mathbb{D}}$, $\alpha = f'(p) \in (0,1)$, then for all $0<\varepsilon<1$, we have $$(f^*)^n(\zeta) \in A(p;\alpha^{(1+\varepsilon)n} , \alpha^{(1-\varepsilon)n}) \cap \partial\mathbb{D}$$ for $n$ large enough and $\lambda$-almost every $\zeta \in \partial\mathbb{D}$.

The strategy  to prove Theorem \ref{A} is as follows. First, we use the functions $M_n$ defined in Lemma \ref{Mn} to convert our autonomous system into a non-autonomous one fixing 0, so that the shrinking targets $D(p,\alpha^{(1-\varepsilon)n}) \cap \partial\mathbb{D}$ and $D(p, \alpha^{(1+\varepsilon)n})\cap\partial\mathbb{D}$ are converted into another sequence of targets of $\partial\mathbb{D}$. As in the proof of \cite[Theorem~E]{benini_shrinking_2024}, the goal will then be to apply Lemma \ref{BEF+A} to these new targets (or to their complements, depending on the case) to prove that our non-autonomous system fails to hit them almost surely. Thus, we need to prove that the new targets are actually shrinking, and that the sum of their lengths is convergent.

\textit{1. Rate of convergence of $f^n(0)$ to $p$}

We are interested in controlling the speed of convergence of $f^n(0)$ to the Denjoy-Wolff point. In this sense, Wolff lemma and the concept of angular derivative are useful tools, since they provide bounds on $|f^n(0) - p|$ in terms of the angular derivative in the hyperbolic case. These bounds are stated in the following lemma. \\

\begin{lemma}\label{cotas+}
    Let $f \colon \mathbb{D} \to \mathbb{D}$ be a hyperbolic self map of $\mathbb{D}$ and let $p \in \partial\mathbb{D}$ be its Denjoy-Wolff point. Let $\alpha = f'(p)$. Then, for every $\delta > 0$, there exists a real constant $C \geq 1$ such that $$\frac{1}{C} (\alpha - \delta)^n \leq |f^n(0) - p| \leq C \alpha^n,$$ for $n$ large enough.

\begin{proof}
To obtain an upper bound, we  apply Lemma \ref{cotas}. To this end, consider $\eta > 1$, so that $0 \in H(p,\eta)$. Hence,
\begin{align*}
    |p-f^n(0)|^2 &\leq \eta (f'(p))^n (1-|f^n(0)|^2) = \eta (f'(p))^n(1+|f^n(0)|) (1-|f^n(0)|)\\ &\leq 2\eta (f'(p))^n  (1-|f^n(0)|) \leq2 \eta (f'(p))^n |p-f^n(0)|,
\end{align*} 
where in the last inequality we have used that, by the triangle inequality, $|p-f^n(0)| \geq 1 - |f^n(0)|$. Dividing both sides by $|p-f^n(0)|$ (which is strictly greater than 0) and denoting $f'(p) = \alpha$, we have \begin{equation*}\label{cotaexphyp} 1- |f^n(0)| \leq |f^n(0) - p| \leq C_1 \alpha^{n},\end{equation*} with $C_1 > 0$.

When it comes to a lower bound, we can easily get one using the definition of the angular derivative, as given that $f$ is hyperbolic, the convergence to the Denjoy-Wolff point is non-tangential, so there exists a Stolz Angle $\Delta$ such that $f^n(0) \in \Delta$ for all $n$ large enough.  As $\alpha = f'(p)$, by Definition \ref{angularderivative} we have that
$$
\lim_{z\to p, z\in\Delta} \frac{f(z) - p}{z-p} = \alpha.
$$ Therefore, for each $\delta >0$ there exists $r > 0$ such that
$$
\alpha - \delta \leq \frac{|f(z)-p|}{|z-p|} \leq \alpha + \delta
$$
for $z \in D(p,r) \cap \Delta$. Moreover, there exists $n_0 \in \mathbb{N}$ such that $f^n(0) \in \Delta \cap D(p,r)$ for all $n \geq n_0$. Thus, for all $n \geq n_0$, $|f^n(0) - p| \geq (\alpha - \delta)^{n-n_0} |f^n(0) - p|  $, which implies that \begin{equation*}\label{lowboundhyp}
    |f^n(0) - p| \geq C_2 (\alpha-\delta)^n,
\end{equation*} with $C_2 > 0$ and for $n \geq n_0$. Taking $C \coloneqq \max\{C_1, \frac{1}{C_2},1\}$, the result follows. 
\end{proof}
\end{lemma}

\textit{2. Two technical lemmas to ensure that the targets shrink}

We now prove two technical lemmas that are necessary to assure that the targets we defined are actually shrinking.\\

\begin{lemma}\label{lemmahyp}
Let $f\colon \mathbb{D} \to \mathbb{D}$ be a hyperbolic function with Denjoy-Wolff point $p\in \partial\mathbb{D}$. For each $n \in \mathbb{N}$, set $r_n > 0$ and let $J_n = D(p,r_n) \cap \partial\mathbb{D}$ be a sequence of shrinking targets verifying

    $$
    \frac{|p-f^n(0)|}{r_n} \to 0, \quad \text{\ as\ }n \to \infty .
    $$

    Then the sequence of targets $I_n = M_n^{-1}(J_n^c)$ shrinks, where $M_n$ are the functions defined in Lemma \ref{Mn}.
 \begin{proof}
Let $J_n = \bigl(e^{i\theta_n},\,e^{i\varphi_n}\bigr) \subset \partial\mathbb{D}$,
and define $\delta_n := e^{i\theta_n} - p$, $
\eta_n := e^{i\varphi_n} - p$, $\gamma_n := f^n(0) - p$.
By assumption, $|\delta_n| =|\eta_n| = r_n$ and $|\gamma_n| = o(r_n)$.

Consider the inverse of the functions $M_n$ defined in Lemma \ref{Mn}, namely
$$
M_n^{-1}(w) = \frac{w - f^n(0)}{1 - \overline{f^n(0)}\,w}
= \frac{w - (p + \gamma_n)}{1 - (\overline{p} + \overline{\gamma_n})\,w}.
$$

Evaluating $M_n^{-1}$ at \(w = e^{i\theta_n} = p + \delta_n\), and taking into account that $\overline{p} p = 1$, we have
\begin{align*}
M_n^{-1}(e^{i\theta_n})
&= \frac{(p + \delta_n) - (p + \gamma_n)}{1 - (\overline{p} + \overline{\gamma_n})(p + \delta_n)} = \frac{\delta_n - \gamma_n}{1 - \overline{p}p - \overline{p}\delta_n - \overline{\gamma_n}p - \overline{\gamma_n}\delta_n} \\ &= \frac{\delta_n - \gamma_n}{- \overline{p}\delta_n - \overline{\gamma_n}p - \overline{\gamma_n}\delta_n} = \frac{\delta_n \left(1-\frac{\gamma_n}{\delta_n}\right)}{-\overline{p}\delta_n \left( 1 
+ \frac{\overline{\gamma_n}p + \overline{\gamma_n}\delta_n}{\overline{p}\delta_n}\right)}.
\end{align*}

Therefore, as $|\gamma_n| = o(r_n) $ and 
$$
\left| \frac{\overline{\gamma_n}p + \overline{\gamma_n}\delta_n }{\overline{p} \delta_n} \right| = \frac{|\gamma_n| |p + \delta_n|}{|\delta_n|} \leq \frac{|\gamma_n| (|p| + |\delta_n|)}{|\delta_n|} = \frac{|\gamma_n|}{|\delta_n|} (1+ |\delta_n|) \to 0, \text{ as }n \to \infty,
$$

it follows that
\[
M_n^{-1}(e^{i\theta_n}) = \frac{\delta_n (1 + o(1))}{-\overline{p} \delta_n (1 + o(1))} = -\frac{1 + o(1)}{\overline{p}} = -(1 + o(1)) p.
\]

Hence
\[
\left| M_n^{-1}(e^{i\theta_n}) + p \right| \to 0.
\]

An analogous argument for \(w = e^{i\varphi_n}\) gives
\[
\left| M_n^{-1}(e^{i\varphi_n}) + p \right| \to 0.
\]

Thus, both endpoints of the arc
\[
I_n^c = M_n^{-1}(J_n) = (e^{i\alpha_n}, e^{i\beta_n})
\]
converge to \(-p\). Since $f^n(0) \to p$ non-tangentially, we have that $M_n^{-1}(p) =  \frac{\gamma_n}{\overline{\gamma_n}p}\not \to -p$ as $n \to \infty$.  
Therefore $|I_n^c| \to 2\pi$, and thus $|I_n| \to 0$.\end{proof}
\end{lemma}\

\begin{lemma}\label{lemmahyp2}
Let $f\colon \mathbb{D} \to \mathbb{D}$ be a hyperbolic function with Denjoy-Wolff point $p$. For each $n \in \mathbb{N}$, set $r_n > 0$ and let $J_n = D(p,r_n) \cap \partial\mathbb{D}$ be a sequence of shrinking targets verifying

    $$
    \frac{r_n}{|p-f^n(0)|} \to 0 \text{\ as\ }n \to \infty .
    $$

    Then the sequence of targets $I_n = M_n^{-1}(J_n)$ shrinks, where $M_n$ are the functions defined in Lemma \ref{Mn}.
 \begin{proof}
Let $J_n = \bigl(e^{i\theta_n},\,e^{i\varphi_n}\bigr) \subset \partial\mathbb{D}$,
and define $\delta_n := e^{i\theta_n} - p$, $
\eta_n := e^{i\varphi_n} - p$, $\gamma_n := f^n(0) - p$.
By assumption, $|\delta_n| = |\eta_n| = r_n$, $|\gamma_n| = o(r_n)$.

Consider the inverse of the functions $M_n$ defined in Lemma \ref{Mn}, namely
$$
M_n^{-1}(w) = \frac{w - f^n(0)}{1 - \overline{f^n(0)}\,w}
= \frac{w - (p + \gamma_n)}{1 - (\overline{p} + \overline{\gamma_n})\,w}.
$$

Evaluating $M_n^{-1}$ at \(w = e^{i\theta_n} = p + \delta_n\), and taking into account that $\overline{p} p = 1$, we have
\begin{align*}
M_n^{-1}(e^{i\theta_n})
&= \frac{(p + \delta_n) - (p + \gamma_n)}{1 - (\overline{p} + \overline{\gamma_n})(p + \delta_n)} = \frac{\delta_n - \gamma_n}{1 - \overline{p}p - \overline{p}\delta_n - \overline{\gamma_n}p - \overline{\gamma_n}\delta_n} \\ &= \frac{\delta_n - \gamma_n}{- \overline{p}\delta_n - \overline{\gamma_n}p - \overline{\gamma_n}\delta_n} = \frac{\gamma_n \left(\frac{\delta_n}{\gamma_n}-1\right)}{-p\overline{\gamma_n} \left( 1 
+ \frac{\overline{p}\delta_n+ \overline{\gamma_n}\delta_n}{p\overline{\gamma_n}}\right)}.
\end{align*}

Therefore, as $|r_n| = |\delta_n| =  o(|\gamma_n|)$ and 
$$
\left| \frac{\overline{p}\delta_n+ \overline{\gamma_n}\delta_n}{p\overline{\gamma_n}} \right| = \frac{|\delta_n| |\overline{p} + \overline{\gamma_n}|}{|\overline{\gamma_n}|} \leq \frac{|\delta_n| (|\overline{p}| + |\overline{\gamma_n}|)}{|\gamma_n|} = \frac{|\delta_n|}{|\gamma_n|} (1+ |\gamma_n|) \to 0, \text{ as }n \to \infty,
$$

it follows that
\[
M_n^{-1}(e^{i\theta_n}) = \frac{\gamma_n (-1 + o(1))}{-p \overline{\gamma_n} (1 + o(1))} = \frac{\gamma_n}{p \overline{\gamma_n}} (1+o(1)).
\]

Hence
\[
\left| M_n^{-1}(e^{i\theta_n}) - \frac{\gamma_n}{p \overline{\gamma_n}} \right| \to 0.
\]

An analogous argument for \(w = e^{i\varphi_n}\) gives
\[
\left| M_n^{-1}(e^{i\varphi_n}) - \frac{\gamma_n}{p \overline{\gamma_n}} \right| \to 0.
\]

Thus, both endpoints of the arc
\[
I_n = M_n^{-1}(J_n) = (e^{i\alpha_n}, e^{i\beta_n})
\]
converge to the moving point \(\frac{\gamma_n}{p \overline{\gamma_n}}\). But $M_n^{-1}(p) = \frac{\gamma_n}{\overline{\gamma_n}p} \in I_n$ as well. 
Therefore $|I_n| \to 0$.\end{proof}
\end{lemma}\

\textit{3. Upper bound on the boundary convergence almost surely}

    We now prove that given $0 < \varepsilon < 1$, $$(f^*)^n(\zeta) \in D(p, \alpha^{(1-\varepsilon)n}) \cap \partial\mathbb{D}$$ for $n$ large enough and $\lambda$-almost every $\zeta \in \partial\mathbb{D}$, where $\alpha = f'(p) \in (0,1)$.
    
    Consider the arc $J_n(\varepsilon) = D(p,\alpha^{(1-\varepsilon)n}) \cap \partial\mathbb{D}$, $n \in \mathbb{N}$, and define 
    \begin{align*}
         E(\varepsilon) &\coloneqq \{ \zeta \in \partial\mathbb{D}\colon f^n(\zeta) \in J_n(\varepsilon) \text{\ for all $n$ large enough}\}\\
    &= \{ \zeta \in \partial\mathbb{D} \colon (f^n(\zeta)) \text{\ fails to hit\ } (J_n(\varepsilon)^c)\},
    \end{align*}
    where $J_n(\varepsilon)^c = \partial\mathbb{D}\setminus J_n(\varepsilon)$. We will prove that $E(\varepsilon)$ has full measure.
    
    We define $G_n = M_n^{-1} \circ f^n$, for $n \in \mathbb{N}$, where $M_n$ are the functions defined in Lemma \ref{Mn}. We can observe that the functions $G_n$ are self maps of $\mathbb{D}$ that fix the origin for all $n$.

    Therefore
    $$
    E(\varepsilon) = \{ \zeta \in \partial\mathbb{D} \colon G_n(\zeta) \text{\ fails to hit\ } (I_n) \},
    $$
    where $I_n = M_n^{-1}(J_n(\varepsilon)^c)$, $n \in \mathbb{N}$.

    We will prove that we can apply Theorem \ref{BEF+A} to the functions $G_n$. By proving this, the result follows immediately, as this would automatically imply that $E(\varepsilon)$ is a set of full measure. Thus, our goal is to see that
    $$
    \sum_{n=1}^\infty |I_n| < \infty.
    $$
    As a result of Lemma \ref{Mn}, we have
    \begin{equation*}\label{senoA}
               2 \sin\left(\frac{1}{2}|I_n|\right) = \frac{(1-|f^n(0)|^2)|e^{i\varphi_n}-e^{i\theta_n}|}{|f^n(0) - e^{i\varphi_n}| |f^n(0) - e^{i\theta_n}|},
    \end{equation*}
    where $J_n(\varepsilon)^c$ is the arc $(e^{i\theta_n},e^{i\varphi_n})$. In order to get an upper bound for this expression, note that we have by Lemma \ref{cotas+} that $1~-~|f^n(0)| \leq |f^n(0) - p| \leq C \alpha^{n}$, with $C > 0$.

    As
    $$
    |f^n(0) - e^{i\varphi_n}| \geq ||f^n(0)-p| - |p-e^{i\varphi}|| \geq |p - e^{i\varphi}| - |f^n(0) - p|,
    $$
    then
    $$
    |f^n(0) - e^{i\varphi_n}| \geq \alpha^{(1-\varepsilon)n} - C\alpha^n,
    $$
    and similarly
    $$
    |f^n(0) - e^{i\theta_n}| \geq \alpha^{(1-\varepsilon)n} - C\alpha^n.
    $$

    Hence, using that $|e^{i\varphi_n} - e^{i\theta_n}| \leq 2 \alpha^{(1-\varepsilon)n}$, that $1+|f^n(0)| <2$ and the above mentioned inequalities, we obtain
    \begin{align*}
        0 &\leq \frac{(1-|f^n(0)|^2)|e^{i\varphi_n} - e^{i\theta_n}|}{|f^n(0) - e^{i\varphi_n}||f^n(0) - e^{i\theta_n}|} =  \frac{(1-|f^n(0)|)(1+|f^n(0)|)|e^{i\varphi_n} - e^{i\theta_n}|}{|f^n(0) - e^{i\varphi_n}||f^n(0) - e^{i\theta_n}|} 
        \\ 
        &\leq \frac{4C\alpha^n  \alpha^{(1-\varepsilon)n}}{\alpha^{2(1-\varepsilon)n} - 2C\alpha^{(1-\varepsilon)n}\alpha^n + C^2\alpha^{2n}} = \frac{4C}{\alpha^{-n\varepsilon} - 2C + C^2 \alpha^{n\varepsilon}} \to 0, \text{\ as $n \to \infty$}.
    \end{align*}

    Hence, we can apply Lemma \ref{lemmahyp} to deduce that, for $n$ large enough, $\frac{1}{2}|I_n| \in [0,\frac{\pi}{2}]$ and therefore
    $$
    \sin\left(\frac{1}{2}|I_n|\right) \geq \frac{2}{\pi} \cdot \frac{1}{2} |I_n| = \frac{1}{\pi} |I_n|.
    $$

    As a result, for $n$ large enough it holds that
    $$
    |I_n| \leq \frac{2C\pi}{\alpha^{-n\varepsilon} - 2C + C^2 \alpha^{n\varepsilon}},
    $$

    where $$\sum_{n=1}^\infty  \frac{2C\pi}{\alpha^{-n\varepsilon} - 2C + C^2 \alpha^{n\varepsilon}} < \infty $$
    by the limit comparison test (comparing with $\sum_{n=1}^\infty \alpha^{n\varepsilon} < \infty $). This implies that $\sum_{n=1}^\infty|I_n| < \infty$, as desired.

    \textit{4. Lower bound on the boundary convergence almost surely}

    We finally prove that given $0 < \varepsilon < 1$, $$(f^*)^n(\zeta) \not \in D(p, \alpha^{(1+\varepsilon)n}) \cap \partial\mathbb{D}$$ for $n$ large enough and $\lambda$-almost every $\zeta \in \partial\mathbb{D}$.

    This proof is very similar to the previous one but considering the arc $J_n(\varepsilon) = D(p,\alpha^{(1+\varepsilon)n}) \cap \partial\mathbb{D}$, $n \in \mathbb{N}$, and 
    \begin{align*}
    E(\varepsilon) &\coloneqq \{ \zeta \in \partial\mathbb{D}\colon f^n(\zeta) \not\in J_n(\varepsilon) \text{\ for all $n$ large enough}\}\\ &=  \{\zeta \in \partial\mathbb{D} \colon (f^n(\zeta)) \text{\ fails to hit\ } (J_n(\varepsilon)) \}  \\ &=       \{ \zeta \in \partial\mathbb{D} \colon G_n(\zeta) \text{\ fails to hit\ } (I_n) \},
    \end{align*}

   where, $G_n = M_n^{-1} \circ f^n$ and  $I_n = M_n^{-1}(J_n(\varepsilon))$, for $n \in \mathbb{N}$.
    We will prove that we can apply Theorem \ref{BEF+A} to the functions $G_n$. As previously, by proving this, the result follows immediately. As a result of Lemma \ref{Mn} we have
    $$
       2 \sin\left(\frac{1}{2}|I_n|\right) = \frac{(1-|f^n(0)|^2)|e^{i\varphi_n}-e^{i\theta_n}|}{|f^n(0) - e^{i\varphi_n}| |f^n(0) - e^{i\theta_n}|},
    $$
    where $J_n(\varepsilon)$ is the arc $(e^{i\theta_n},e^{i\varphi_n})$.

    By Lemma \ref{cotas+} we have that $|p-f^n(0)| \geq C\alpha^{(1+\frac{\varepsilon}{3})n}$, with $C > 0$, for $n$ large enough.

    As
    $$
    |f^n(0) - e^{i\varphi_n}| \geq ||f^n(0)-p| - |p-e^{i\varphi}|| \geq  |f^n(0) - p| -|p - e^{i\varphi}|,
    $$
    then
    $$
    |f^n(0) - e^{i\varphi_n}| \geq  C\alpha^{(1+\frac{\varepsilon}{3})n} - \alpha^{(1+\varepsilon)n},
    $$
    and similarly
    $$
    |f^n(0) - e^{i\theta_n}| \geq C\alpha^{(1+\frac{\varepsilon}{3})n} - \alpha^{(1+\varepsilon)n}.
    $$

    Hence, using that $|e^{i\varphi_n} - e^{i\theta_n}| \leq 2 \alpha^{(1+\varepsilon)n}$, that $1+|f^n(0)| <2$ and the above mentioned inequalities, we obtain
    \begin{align*}
        0 &\leq \frac{(1-|f^n(0)|^2)|e^{i\varphi_n} - e^{i\theta_n}|}{|f^n(0) - e^{i\varphi_n}||f^n(0) - e^{i\theta_n}|} =  \frac{(1-|f^n(0)|)(1+|f^n(0)|)|e^{i\varphi_n} - e^{i\theta_n}|}{|f^n(0) - e^{i\varphi_n}||f^n(0) - e^{i\theta_n}|} 
        \\  
        &\leq \frac{4C\alpha^n  \alpha^{(1+\varepsilon)n}}{\alpha^{2(1+\varepsilon)n} - 2C\alpha^{(1+\varepsilon)n}\alpha^{(1+\frac{\varepsilon}{3})n} + C^2\alpha^{2(1+\frac{\varepsilon}{3})n}} = \frac{4C}{\alpha^{n\varepsilon} - 2C\alpha^{\frac{\varepsilon}{3}\cdot n} + C^2 \alpha^{-\frac{1}{3} \cdot \varepsilon n}} \to 0, 
    \end{align*}
 as $n \to \infty$.
 
    Then we can apply Lemma \ref{lemmahyp2} to deduce that, for $n$ large enough, $\frac{1}{2}|I_n| \in [0,\frac{\pi}{2}]$ and therefore
    $$
    \sin\left(\frac{1}{2}|I_n|\right) \geq \frac{2}{\pi} \cdot \frac{1}{2} |I_n| = \frac{1}{\pi} |I_n|.
    $$

    As a result, for $n$ large enough it holds that
    $$
    |I_n| \leq \frac{2C \pi}{\alpha^{n\varepsilon} - 2C\alpha^{\frac{\varepsilon}{3}\cdot n} + C^2 \alpha^{-\frac{1}{3} \cdot \varepsilon n}},
    $$

    where $$\sum_{n=1}^\infty  \frac{2C\pi}{\alpha^{n\varepsilon} - 2C\alpha^{\frac{\varepsilon}{3}\cdot n} + C^2 \alpha^{-\frac{1}{3} \cdot \varepsilon n}} < \infty $$
    by the limit comparison test (comparing with $\sum_{n=1}^\infty \alpha^{\frac{1}{3}\cdot n\varepsilon} < \infty $). This implies that $\sum_{n=1}^\infty|I_n| < \infty$. The proof is complete. \hfill $\square$

\bibliography{refs}
\end{document}